\documentclass[amstex,12pt,english,amssymb]{article}

\usepackage{mathtext}
\usepackage[cp1251]{inputenc}
\usepackage[T2A]{fontenc}
\usepackage[english]{babel}
\usepackage[dvips]{graphicx}
\usepackage{amsmath}
\usepackage{amssymb}
\usepackage{amsxtra}
\usepackage{latexsym}
\usepackage{ifthen}

\usepackage[dvips]{graphicx}
\usepackage{epstopdf}
\usepackage{epsfig}
\usepackage{epic}
\usepackage{eepic}

\newenvironment{proof}[1][Proof]{\emph{#1.} }{\ \rule{0.5em}{0.5em}}

\begin{document}




\numberwithin{equation}{section}
\newcommand{\abs}[1]{\lvert#1\rvert}
\newcommand{\blankbox}[2]{%
  \parbox{\columnwidth}{\centering
    \setlength{\fboxsep}{0pt}%
    \fbox{\raisebox{0pt}[#2]{\hspace{#1}}}
} }

\newcounter{lemma}[section]
\newcommand{\lemma}{\par \refstepcounter{lemma}%
{\bf Lemma \arabic{section}.\arabic{lemma}.}}
\renewcommand{\thelemma}{\thesection.\arabic{lemma}}

\newcounter{corollary}[section]
\newcommand{\corollary}{\par \refstepcounter{corollary}%
{\bf Corollary \arabic{section}.\arabic{corollary}.}}
\renewcommand{\thecorollary}{\thesection.\arabic{corollary}}

\newcounter{remark}[section]
\newcommand{\remark}{\par \refstepcounter{remark}%
{\bf Remark \arabic{section}.\arabic{remark}.}}
\renewcommand{\theremark}{\thesection.\arabic{remark}}

\newcounter{theorem}[section]
\newcommand{\theorem}{\par \refstepcounter{theorem}%
{\bf Theorem \arabic{section}.\arabic{theorem}.}}
\renewcommand{\thetheorem}{\thesection.\arabic{theorem}}

\newcounter{proposition}[section]
\newcommand{\proposition}{\par \refstepcounter{proposition}%
{\bf Proposition \arabic{section}.\arabic{proposition}.}}
\renewcommand{\theproposition}{\thesection.\arabic{proposition}}

\def\N{{{\Bbb N}}}
\def\Z{{{\Bbb Z}}}
\def\T{{{\Bbb T}}}
\def\R{{\Bbb R}}
\def\D{\Delta}
\def\Q{{\Bbb Q}}
\def\C{{\Bbb C}}
\def\l{{\lambda }}
\def\a{{\alpha }}
\def\f{{\theta}}
\def\r{{\rho}}
\def\a{{\alpha}}
\def\b{{\beta}}
\def\d{{\delta}}
\def\e{{\varepsilon}}
\def\s{{\sigma}}
\def\supp{\mbox{\rm supp}}

\def\Xint#1{\mathchoice
   {\XXint\displaystyle\textstyle{#1}}%
   {\XXint\textstyle\scriptstyle{#1}}%
   {\XXint\scriptstyle\scriptscriptstyle{#1}}%
   {\XXint\scriptscriptstyle\scriptscriptstyle{#1}}%
   \!\int}
\def\XXint#1#2#3{{\setbox0=\hbox{$#1{#2#3}{\int}$}
     \vcenter{\hbox{$#2#3$}}\kern-.5\wd0}}
\def\dashint{\Xint-}

\def\cc{\setcounter{equation}{0}
\setcounter{figure}{0}\setcounter{table}{0}}

\overfullrule=0pt

\title{{\bf On approximation methods generated by\\ generalized Bochner-Riesz kernels}}

\author{{\bf Yu.\,S.~Kolomoitsev}\\}
\date{\today}
\maketitle

\begin{abstract}
Some sharp results related to the convergence of means and families
of operators generated by the generalized Bochner-Riesz kernels are
obtained. The exact order of approximation of functions by these
methods via $K$-functional (or its realization in the case of the
space $L_p$, $0<p<1$) is derived.
\end{abstract}

\large

\section{Introduction} Let $\T^d=[0,2\pi)^d$ be the $d$-dimensional torus. As usual, the space
$L_p(\T^d)$, $0<p<\infty$, consists of measurable real valued
functions $f(x)$, $x\in\R^d$, which is $2\pi$-periodic in each
variable and
$$
\Vert f\Vert_p=\bigg(\int_{\T^d}|f(x)|^p dx\bigg)^{\frac 1p}<\infty.
$$
By $L_\infty(\T^d)$ denote the space of all real valued
$2\pi$-periodic continuous functions on $\T^d$ which is equipped
with the norm
$$
\Vert f\Vert_\infty=\max_{x\in\T^d}|f(x)|.
$$

Binomial coefficients of order $\b>0$ are given by
$$
\binom{\b}{k}=\frac{\b(\b-1)\dots(\b-k+1)}{k!},\quad k\in\N,
$$
and $\binom{\b}{0}=1$. We will repeatedly use the following
estimates (see, for example,~\cite[Ch.1,\,\S 1]{SKM})
\begin{equation}\label{eqBinom}
\bigg|\binom{\b}{k}\bigg|\le \frac{C(\b)}{k^{\b+1}}.
\end{equation}

We denote by $C$ and $C_j$, $j=1,2,\dots,$  positive constants
depending on the indicated parameters.  The notation $A(f, n) \asymp
B(f, n) $ means  a two-sided inequality with positive constants
independent of $ f $ and $ n $.

The generalized Bochner-Riesz kernel is defined as follows
\begin{equation}\label{eqBRK1}
\mathcal{R}_n^{\b,\,\d}(x)=\sum_{k\in
\Z^d}\bigg(1-\frac{|k|^\b}{n^\b}\bigg)_+^\d e^{i(k,\,x)},
\end{equation}
where $(k,x)=k_1x_1+\dots+k_dx_d$, $|x|=(x,x)^{\frac 12}$, and
$x_+=\max\{x,0\}$. If $\b=2$, then in \eqref{eqBRK1} we have the
classical Bochner-Riesz kernel.

%

In the present paper we deal with the generalized Bochner-Riesz
means given by
\begin{equation}\label{eqFMBR}
\mathcal{S}_n^{\b,\,\d}(f;x)=(2\pi)^{-d}\int_{\T^d}f(x+y)\mathcal{R}_n^{\b,\,\d}(y)dy,\quad
n\in\N,
\end{equation}
and with the family of linear polynomials operators given by
\begin{equation}\label{eqLPBR}
\mathcal{S}_{n;\,\l}^{\b,\,\d}(f;x)=(2n+1)^{-d}\sum_{k=0}^{2n}f(t_n^k+\l)\mathcal{R}_n^{\b,\,\d}(x-t_n^k-\l),\quad
n\in\N,
\end{equation}
where
$$
t_n^k=\frac{2\pi k}{2n+1},\quad k\in\Z^d;\quad
\sum_{k=0}^{2n}=\sum_{k_1=0}^{2n}\dots\sum_{k_d=0}^{2n}.
$$
The Bochner-Riesz means and the family of polynomial operators
generated by the kernel $\mathcal{R}_n^{2,\,\d}$ we denote by
$\mathcal{S}_n^{\d}$ and $\{\mathcal{S}_{n;\,\l}^{\d}\}$,
respectively.


The properties of the classical Bochner-Riesz means
$\mathcal{S}_n^{\d}$ are intensively studied by many authors (see,
for example,~\cite[Ch.7]{StWE}, \cite[Ch.3]{Graf}, \cite{Gol82},
\cite{Bel75}, \cite{KuzTr80}, \cite{Tr80}). Some sharp results
related to the approximation of functions by these means
$\mathcal{S}_n^{\d}$ and by the family
$\{\mathcal{S}_{n;\,\l}^{\d}\}$ were obtained in~\cite{RS08}.

Which approximation properties do the means~\eqref{eqFMBR} have for
different values of  $\b$ and $\d$? It is known that the regularity
(convergence) of the Bochner-Riesz means depends on the parameter
$\d$. In particular, it was shown in~\cite{Tr80} and~\cite{KuzTr80}
that if $\d>{(d-1)}/{2}$,  then the means~\eqref{eqFMBR} converge in
$L_p(\T^d)$ for any $p\in [1,\infty]$, for other values of $\d$ the
convergence may not be achieved. The parameter $\b>0$ does not have
any effect on the regularity of the means $\mathcal{S}_n^{\b,\,\d}$
in contrast to $\d$ (see, for example,~\cite[Ch.8]{TrigBook} in the
case of even $\b$). However, the order of approximation of $f$ by
the means $\mathcal{S}_n^{\b,\,\d}(f)$ is better for large $\b$ and
the order is the same while $\d$ is changing.

The special module of smoothness of order $\b>0$ is defined by
\begin{equation*}
\widetilde{\omega}_\b(f,h)_\infty=\bigg\Vert \int_{|u|\ge
1}\sum_{\nu=0}^{2r}\binom{2r}{\nu}(-1)^\nu f(\cdot +
(\nu-r)uh)\frac{du}{|u|^{d+\b}}\bigg\Vert_\infty,
\end{equation*}
where $r\in\N$ and $r>d-1+\b$. It was shown in~\cite{KuzTr80} that
the approximation error of $f\in C(\T^d)$ by
$\mathcal{S}_n^{\b,\,\d}(f)$ is equivalent to
$\widetilde{\omega}_\b(f,h)_\infty$ under certain restrictions on
$\b$ and $\d$. This result is also valued in the space $L_p(\T^d)$
(see, for example,~\cite[Theorem 7]{Tr80}). In particular, we have
the following theorem.

\medskip

\noindent{\textbf{Theorem A.}} \emph{Let $f\in L_p(\T^d)$, $1\le
p\le\infty$, $\b>0$, and $\d>{(d-1)}/{2}$. Then
\begin{equation}\label{eqA.0}
\Vert f-\mathcal{S}_{n}^{\b,\,\d}(f)\Vert_{{p}}\asymp
\widetilde{\omega}_\b(f,1/n)_p\,,\quad n\in\N.
\end{equation}}


Note that in the case of even $\b>0$ in~\eqref{eqA.0} the module of
smoothness $\widetilde{\omega}_\b(f,h)_p$ can be replaced by the
corresponding $K$-functional (see~\cite[Ch.8]{TrigBook}) and other
special moduli of smoothness (see~\cite{Tr80}).

In~\cite{RS08} it was shown that for the classical Bochner-Riesz
means $\mathcal{S}_n^{\d}$ as well as for the families
$\{\mathcal{S}_{n;\,\l}^{\d}\}$ there is an alternative:
\emph{either the means  $\mathcal{S}_n^{\d}$ (or the family
$\{\mathcal{S}_{n;\,\l}^{\d}\}$) diverge in $L_p$ or its
approximation error is equivalent to the $K$-functional (or its
realization if $0<p<1$).}

In the present paper the results of~\cite{RS08} are extended to the
case of the generalized Bochner-Riesz means
$\mathcal{S}_n^{\b,\,\d}$ as well as to the family
$\{\mathcal{S}_{n;\,\l}^{\b,\,\d}\}$ with any $\b\in\R$, $\b>0$. In
particular, it is proved that the above alternative holds for any
positive $\d$ and $\b$. It turns out that in the case of $0<p<1$ the
regularity of the family $\{\mathcal{S}_{n;\,\l}^{\b,\,\d}\}$
essentially depends on the parameter $\b$.

The paper is organized as follows. In Section~\ref{s2} we formulate
the main results. In Subsection~\ref{s2.1} the convergence theorems
are formulated; in Subsection~\ref{s2.2} the theorems on equivalence
of error of approximation of functions by corresponding method are
formulated.  The auxiliary results are formulated and proved in
Section~\ref{s3}. In Section~\ref{s4} we prove the main results of
the paper.

\section{Main Results}\label{s2}

\subsection{Convergence Theorems}\label{s2.1}
Often we will deal with functions in $L_p(\T^{2d})$ which depend
additionally on parameter $\lambda\in\T^d$.
We denote by $\Vert\cdot\Vert_{\overline{p}}$  the $p$-(quasi-)norm
with respect to both the main variable $x\in\T^d$ and the parameter
$\lambda\in\T^d$, i.e.
$$
\Vert\cdot\Vert_{\overline{p}}=\Vert
\Vert\cdot\Vert_{p;\,x}\Vert_{p;\,\l}\,,
$$
where $\Vert\cdot\Vert_{p;\,x}$ and $\Vert\cdot\Vert_{p;\,\l}$ are
the $p$-norms (quasi-norms, if $0<p<1$) with respect to $x$ and
$\l$, respectively.

Let $\mathcal{T}_n$ be the set of all real valued trigonometric
polynomials of order $n$:
\begin{equation*}
\mathcal{T}_n=\bigg\{T(x)=\sum_{k\in\Z^d,~|k|\le n} c_ke^{i(k,x)}~:~
c_{-k}=\overline{c}_k\bigg\}.
\end{equation*}
A sequence of linear operators $\{\mathcal{L}_n\}_{n\in\N}$, mapping
$L_p$, $1\le p\le\infty$, into the space $\mathcal{T}_n$ is said to
be convergent (or converges) in $L_p$, if
\begin{equation*}
\lim_{n\to\infty}\Vert f-\mathcal{L}_n(f)\Vert_p=0
\end{equation*}
for each $f\in L_p(\T^d)$. By analogy, a family of linear operators
$\{\mathcal{L}_{n;\,\lambda}\}_{n\in\N,\, \lambda\in\T^d}$, mapping
$L_p$, $0< p\le\infty$, into $\mathcal{T}_n$, converges in $L_p$, if
\begin{equation*}
\lim_{n\to\infty}\Vert
f-\mathcal{L}_{n;\,\lambda}(f)\Vert_{\overline{p}}=0
\end{equation*}
for each $f\in L_p(\T^d)$.

In order to formulate the main results, we split the domain $\R_+^2$
of pairs $(1/p,\b)$ into three parts:
\begin{equation*}
\Sigma(d)=\bigg\{\bigg(\frac1p,\d\bigg)\in\R_+^2~:~\d>\max\bigg\{\frac{d-1}2,d\bigg(\frac1p-\frac12\bigg)-\frac12\bigg\}
\bigg\},
\end{equation*}
\begin{equation*}
\Gamma(d)=\bigg\{\bigg(\frac1p,\d\bigg)\in\R_+^2~:~0\le\d\le
d\bigg|\frac1p-\frac12\bigg|-\frac12\bigg\},\qquad\qquad\quad
\end{equation*}
\begin{equation*}
\Omega(d)=\bigg\{\bigg(\frac1p,\d\bigg)\in\R_+^2~:~0\le\d\le\frac{d-1}2,~\d>d\bigg|\frac1p-\frac12\bigg|-\frac12
\bigg\}.
\end{equation*}


\begin{theorem}\label{thSh1}
\emph{Let $1\le p\le \infty$ and $\b,\d>0$.  The means
$\mathcal{S}_n^{\b,\,\d}$ converge in  $L_p$ if and only if the
means $\mathcal{S}_n^{\d}$ converge in the same space $L_p$. In
particular, the means $\mathcal{S}_n^{\b,\,\d}$ converge in $L_p$ if
$(1/p,\d)\in \Sigma(d)$ and they diverge in $L_p$ if $(1/p,\d)\in
\Gamma(d)$.}
\end{theorem}

Note that the problem of convergence in the domain $\Omega(d)$ is
not completely studied even for the classical Bochner-Riesz means
$\mathcal{S}_n^{\d}$ (see, for example,~\cite{St}).

From Theorem~\ref{thSh1} it follows that the parameter $\b>0$ does
not affect on the convergence of the means $\mathcal{S}_n^{\b,\,\d}$
in $L_p$, $1\le p\le\infty$. In contrast to the case $p\ge 1$, the
dependence on the parameter $\b$ is essential for the
family~\eqref{eqLPBR} in $L_p$, $0<p<1$. To state the next theorem
we introduce the following set
$$
\mathrm{B}(d)=\bigg\{\bigg(\frac1p,\b\bigg)\in\R_+^2~:~\b\in2\N,\quad
\b>d\bigg(\frac 1p-1\bigg)_{\!\!+}\bigg\}.
$$

\begin{theorem}\label{thSh2}
\emph{1) Let $1\le p\le \infty$ and $\b,\,\d>0$. Then the family
$\{\mathcal{S}_{n;\,\l}^{\b,\,\d}\}$ converges in $L_p$ if and only
if the family $\{\mathcal{S}_{n;\,\l}^{\d}\}$ converges in the same
space $L_p$. In particular, the family
$\{\mathcal{S}_{n;\,\l}^{\b,\,\d}\}$ converges in $L_p$ if
$(1/p,\d)\in \Sigma(d)$ and it diverges in $L_p$ if $(1/p,\d)\in
\Gamma(d)$. }

\emph{2) Let $0<p<1$ and $\b,\,\d>0$. Then the family
$\{\mathcal{S}_{n;\,\l}^{\b,\,\d}\}$ converges in $L_p$ if and only
if $(1/p,\d)\in\Sigma(d)$ and $(1/p,\b)\in \mathrm{B}(d)$.}
\end{theorem}

\subsection{Two-sided estimates of approximation}\label{s2.2}
For our purpose we will use a $K$-functional related to the power of
the Laplacian  $\Delta^{\b/2}$, which we define by
$$
\D^{\b/2}f(x)\sim\sum_k |k|^{\b}c_k(f)e^{-i(k,x)},\quad
\b\in\R,\quad \b>0,
$$
where
$$
c_k(f)=(2\pi)^{-d}\int_{\T^d}f(x)e^{-i(k,x)}dx,\quad k\in\Z^d,
$$
are the Fourier coefficients of the function $f$. The corresponding
$K$-functional is given by
\begin{equation}\label{eqKf1}
K_{\b}(f,t)_p=\inf_g\{\Vert f-g\Vert_p+t^{\b}\Vert\D^{\b/2}
g\Vert_p\}.
\end{equation}

It should be noticed that in the case $0<p<1$ the $K$-functional
given in~\eqref{eqKf1} is identically zero (see \cite{DHI}).
However, in accordance with the concept of ``Realization''
(see~\cite{HI90}), the $K$-functional can be replaced by the
quantity
\begin{equation}\label{eqKf2}
\widetilde{K}_{\b}(f,t)_p=\inf_{T\in\mathcal{T}_{1/t}}\{\Vert
f-T\Vert_p+t^{\b}\Vert \D^{\b/2} T\Vert_p\}
\end{equation}
(see also~\cite{DHI} for properties of $\widetilde{K}_{\b}$ in the
case $0<p<1$).


\begin{theorem}\label{thKf1}
\emph{Let $f\in L_p(\T^d)$, $1\le p\le\infty$, $\b>0$, and
$\d>{(d-1)}/{2}$. Then
\begin{equation}\label{eq.thKf1.0}
\begin{split}
\Vert f-\mathcal{S}_{n;\,\l}^{\b,\,\d}(f)\Vert_{\overline{p}}&\asymp
\Vert f-\mathcal{S}_{n}^{\b,\,\d}(f)\Vert_{{p}}\asymp
K_\b(f,1/n)_p\\
&\asymp\widetilde{K}_\b(f,1/n)_p\asymp
\widetilde{\omega}_\b(f,1/n)_p\,.
\end{split}
\end{equation}}
\end{theorem}

The next theorem holds in the case $0<p<1$.

\begin{theorem}\label{thKf2}
\emph{Let $f\in L_p(\T^d)$, $0<p<1$, $(1/p,\d)\in\Sigma(d)$, and
$(1/p,\b)\in \mathrm{B}(d)$. Then
\begin{equation}\label{eq.thKf2.0}
\Vert f-\mathcal{S}_{n;\,\l}^{\b,\,\d}(f)\Vert_{\overline{p}}\asymp
\widetilde{K}_\b(f,1/n)_p\,,\quad n\in\N.
\end{equation}}
\end{theorem}

The following theorems contain approximation properties of the
methods \eqref{eqFMBR} and \eqref{eqLPBR} for
$(1/p,\d)\in\Omega(d)$. We will denote the subdomains of
$\Omega(d)$, where the Bochner-Riesz means $\mathcal{S}_{n}^{\d}$
and the corresponding families $\{\mathcal{S}_{n;\,\l}^{\d}\}$ have
the convergence property with $\Omega'(d)$  and $\Omega''(d)$,
respectively.

\begin{theorem}\label{thOm1}
\emph{Let $f\in L_p(\T^d)$, $1\le p\le \infty$, $\b,\,\d>0$, and
$(1/p,\d)\in\Omega(d)$. Then
\begin{equation}\label{eqOm1}
\Vert f-\mathcal{S}_n^{\b,\,\d}(f)\Vert_p\le (2\pi)^{-\frac dp}\Vert
f-\mathcal{S}_{n;\,\l}^{\b,\,\d}(f)\Vert_{\overline{p}}\,,\quad
n\in\N,
\end{equation}
in particular, $\Omega''(d)\subset\Omega'(d)$.}
\end{theorem}

\begin{theorem}\label{thOm2}
\emph{Let $f\in L_p(\T^d)$, $1\le p\le \infty$, $d>1$, $\b,\,\d>0$,
and $(1/p,\d)\in\Omega'(d)$. Then
\begin{equation}\label{eqOm2}
\begin{split}
\Vert f-\mathcal{S}_{n}^{\b,\,\d}(f)\Vert_{{p}}&\asymp K_\b(f,1/n)_p
\asymp\widetilde{K}_\b(f,1/n)_p \\
&\asymp \widetilde{\omega}_\b(f,1/n)_p\,.
\end{split}
\end{equation}}
\end{theorem}

Thus, from Theorems~\ref{thSh1} and~\ref{thOm2} it follows that
\emph{the approximation error of the means $\mathcal{S}_n^{\b,\,\d}$
is equivalent to the corresponding $K$-functional if and only if the
means $\mathcal{S}_n^{\b,\,\d}$ converge in $L_p$.}

\begin{theorem}\label{thOm3}
\emph{Let $f\in L_p(\T^d)$, $1\le p\le \infty$, $d>1$, $\b,\,\d>0$,
and $(1/p,\d)\in\Omega''(d)$. Then
\begin{equation}\label{eqOm3}
\begin{split}
\Vert f-\mathcal{S}_{n;\,\l}^{\b,\,\d}(f)\Vert_{\overline{p}}&\asymp
\Vert f-\mathcal{S}_{n}^{\b,\,\d}(f)\Vert_{{p}}\asymp
K_\b(f,1/n)_p\\
&\asymp\widetilde{K}_\b(f,1/n)_p\asymp
\widetilde{\omega}_\b(f,1/n)_p\,.
\end{split}
\end{equation}}
\end{theorem}

Similarly to the case considered above, from Theorems~\ref{thSh2}
and~\ref{thOm3} it follows that \emph{the approximation error of the
family $\{\mathcal{S}_{n;\,\l}^{\b,\,\d}\}$ is equivalent to the
corresponding $K$-functional (or its realization $\widetilde{K}$ if
$0<p<1$) if and only if the family
$\{\mathcal{S}_{n;\,\l}^{\b,\,\d}\}$ converges in $L_p$.}

\section{Auxiliary assertions}\label{s3}


Let us present some facts related to multipliers for trigonometric
polynomials. Let $g$ be a real or complex valued function defined on
$\R^d$. It generates operators $\{A_n(g)\}_{n\ge 1}$ given by
\begin{equation*}\label{eqSP}
A_n(g)T(x)=\sum_{k\in\Z^d}g\bigg(\frac k n\bigg)c_ke^{i(k,x)},\quad
T(x)=\sum_{k\in\Z^d}c_ke^{i(k,x)}\in\mathcal{T},
\end{equation*}
where $\mathcal{T}$ is the set of all trigonometric polynomials.

Consider the inequality
\begin{equation}\label{eqNer1}
\Vert A_n(g)T\Vert_p\le C\Vert T\Vert_p\,,\quad
T\in\mathcal{T}_n,\,\,\,n\ge 1.
\end{equation}
We say that (\ref{eqNer1}) is valid for the function $g$ (this we
denote by $g\in M_p(\mathcal{T})$), if it is valid in the $L_p$-norm
for all $T\in\mathcal{T}_n$ and $n\ge 1$ with some positive constant
$C$ independent of $T$ and $n$.

The following two lemmas are evident.
\begin{lemma}\label{lemSvMp}
\emph{Let $g,h\in M_p(\mathcal{T})$. Then the functions $g+h$ and
$g\cdot h$ belong to $M_p(\mathcal{T})$.}
\end{lemma}

We will also use the inequalities of type
\begin{equation}\label{eqNerr2}
\Vert A_n(g)T\Vert_p\le C\Vert A_n(h)T\Vert_p\,,\quad
T\in\mathcal{T}_n,\,\,\,n\ge 1.
\end{equation}
In the next we suppose that $h(\xi)\neq 0$ for $\xi\neq 0$. Put
$$
\mathcal{X}(\xi)=\frac{g(\xi)}{h(\xi)},\quad
\xi\in\R^d\setminus\{0\}.
$$
We assume that $\mathcal{X}$ is somehow defined at the point
$\xi=0$.

\begin{lemma}\label{lemPSt}
\emph{Let $g(0)=h(0)=0$ and let $\mathcal{X}\in M_p(\mathcal{T})$.
Then the inequality in \eqref{eqNerr2} is valid in $L_p$
independently of the value $\mathcal{X}(0)$.}
\end{lemma}

As usual, the Fourier transform of a function $f\in L_1(\R^d)$ is
given by
$$
\widehat{f}(x)=\int_{\R^d}f(y)e^{-i(x,y)}dy.
$$

The next lemma (see~\cite{Tr80},~\cite[p.\,150-151]{ST}) gives
sufficient conditions for the validity of~\eqref{eqNer1} in the
space $L_p(\T^d)$.

\begin{lemma}\label{lemDostNer}
\emph{Let $0<p\le\infty$ and let  $g$ be a continuous function with
compact support. If $\widehat{g}\in L_{p^*}(\R^d)$
$(p^*=\min(1,p))$, then $g\in M_p(\mathcal{T})$.}
\end{lemma}

Let us denote by $W_1^m(\R^d)$ the Sobolev space of all integrable
functions whose derivatives up to the order  $m$ belong to
$L_1(\R^d)$. By ${\overset{\circ}{W}}\!\!\phantom{l}_1^m(\R^d)$
denote the set of functions in $W_1^m(\R^d)$ having compact support
(see details, for example, in~\cite[Ch.1]{Trie}).

The following lemma is proved in \cite{RS08} (more general
statements are proved in~\cite[Ch.6]{TrigBook}).
\begin{lemma}\label{lemDer}
\emph{Let $0<p\le\infty$, $m=[d/ {p^*}]+1$ $(p^*=\min(1,p)$, $[a]$
is integral part of $a)$. If
$g\in{\overset{\circ}{W}}\!\!\phantom{l}_1^m(\R^d)$, then
$\widehat{g}\in L_{{p^*}}(\R^d)$ and, therefore, $g\in
M_p(\mathcal{T})$.}
\end{lemma}

In the next by $\mathfrak{C}^d$ we denote the class of real or
complex valued $C^\infty$-functions with a compact support contained
in the set $\{x\in\R^d\,:\,|x|\le 1\}$. We use the symbol
$\mathfrak{R}^d$ to denote the class of real valued radial
$C$-functions $\psi$ with a compact support and $\psi(0)=1$.

The following two lemmas are proved in~\cite{RS01}.

\begin{lemma}\label{lemRuOdf1}
\emph{Suppose $f\in C^\infty(\R^d\setminus\{0\})$ is a homogeneous
function of order $\b>0$, it is not a polynomial, and $h\in
\mathfrak{R}^d$. Then $\widehat{fh}\in L_p(\R^d)$ if and only if
$\frac{d}{d+\b}<p\le \infty$.}
\end{lemma}


\begin{lemma}\label{lemRuOdfBer}
\emph{Let $0<p\le \infty$, $\b\in \mathrm{B}(d)$, $n\in \N$, and
$T\in\mathcal{T}_n$. Then
$$
\Vert \Delta^{\b/2}T\Vert_p\le Cn^\b \Vert T\Vert_p\,,
$$
where $C$ is a constant independent of  $T$ and $n$.}
\end{lemma}

The following lemma is easily obtained by using Theorem~6.1.1
in~\cite{TrigBook}.

\begin{lemma}\label{lemProdF}
\emph{Let $0<p\le 1$, let $f$ and $g$ be bounded functions with
compact supports contained in $\{x\in\R^d\,:\,|x|< 3\}$, and let
$\widehat{f}$, $\widehat{g}\in$ $L_p(\R^d)$. Then $\widehat{fg}\in
L_p(\R^d)$.}
\end{lemma}

The item 1) of the next lemma is well known (see, for
example,~\cite{Herz}, \cite[Ch.9]{St}); the item 2) is proved
in~\cite{RS08}.

\begin{lemma}\label{lemBR2}
\emph{1) Let $1\le p\le \infty$. Then the means $\mathcal{S}_n^\d$
converge in $L_p$ if $(1/p,\d)\in \Sigma(d)$ and they diverge in
$L_p$ if $(1/p,\d)\in \Gamma(d)$.}

\emph{2) Let $0<p\le\infty$. Then the family
$\{\mathcal{S}_{n;\,\l}^{\d}\}$ converges in $L_p$ if $(1/p,\d)\in
\Sigma(d)$ and it  diverges in $L_p$ if $(1/p,\d)\in \Gamma(d)$.}
\end{lemma}

Throughout what follows we will use the following notation:
\begin{equation}\label{eqPhi}
\varphi_{\b,\,\d}(x)=(1-|x|^\b)_+^\d,
\end{equation}
\begin{equation}\label{eqh0h1h2}
\begin{split}
&h_0(x)=\left\{
         \begin{array}{ll}
           1, & \hbox{$|x|\le  4/3$;} \\
           0, & \hbox{$|x|>2$,}
         \end{array}
       \right.\quad
h_1(x)=\left\{
         \begin{array}{ll}
           1, & \hbox{$|x|\le 1/2$;} \\
           0, & \hbox{$|x|\ge 3/4$,}
         \end{array}
       \right.\\
&h_2(x)=h_0(x)-h_1(x).
\end{split}
\end{equation}
In addition, suppose that $h_0$ and $h_1$ belong to
$C^\infty(\R^d)\cap \mathfrak{R}^d$.

\begin{lemma}\label{lemPsiBes}
\emph{For $\d>0$ we have
\begin{equation}
\begin{split}
\int_{\R^d}h_2(x)&\varphi_{2,\d}(x) e^{-i(x,y)}dx=\\
&=\sqrt{\frac 2\pi}\cdot\frac{\cos(2\pi|y|-\frac{\pi\d}{2}-\frac\pi
4)}{|y|^{\frac{d+1}{2}+\d}}+O(|y|^{-\frac{d+3}{2}-\d}),\quad |y|\ge
1.
\end{split}
\end{equation}}
\end{lemma}

\begin{proof}
Using the well-known equality (see, for example,~\cite[Ch.
IV]{StWE})
\begin{equation}\label{eqBesrav}
\int_{\R^d}(1-|x|^2)_+^\d
e^{-i(x,y)}dx={\pi^{-\d}}{\Gamma(\d+1)}\frac{J_{\frac d2+\d}(2\pi
|y|)}{|y|^{\frac d2+\d}},
\end{equation}
and $(h_2(x)-1)(1-|x|^2)_+^\d\in C^\infty(\R^d)$, we get
\begin{equation*}
\begin{split}
&\int_{\R^d}h_2(x)\varphi_{2,\d}(x)
e^{-i(x,y)}dx=\\
&=\int_{\R^d}\varphi_{2,\d}(x)
e^{-i(x,y)}dx+\int_{\R^d}(h_2(x)-1)\varphi_{2,\d}(x)
e^{-i(x,y)}dx=\\
&={\pi^{-\d}}{\Gamma(\d+1)}\frac{J_{\frac d2+\d}(2\pi
|y|)}{|y|^{\frac d2+\d}}+O(|y|^{-r}),
\end{split}
\end{equation*}
where $r$ is large enough. It remains to use the asymptotics of the
Bessel function (see, for example,~\cite[Ch. IV]{StWE})
\begin{equation}\label{eqBesAs}
J_\nu(u)=\sqrt{\frac{2}{\pi
u}}\cos\left(u-\frac{\pi\nu}{2}-\frac{\pi}{4}\right)+O(u^{-\frac
32}),\quad u\to\infty.
\end{equation}
\end{proof}

\medskip

The statement of the next lemma for $p=1$ see in~\cite{Lf} or
\cite{Tr77}; for $\b=2$ this result follows immediately from the
asymptotics of Bessel functions \eqref{eqBesAs}.

\begin{lemma}\label{lemFphi}
\emph{Let $0<p\le 1$, $\b>0$, and $\d>0$. Then
${\widehat{\varphi}}_{\b,\,\d}\in L_p(\R^d)$ if and only if
$(1/p,\b)\in \mathrm{B}(d)$ and $\d>d(1/p-1/2)-1/2$.}
\end{lemma}

\begin{proof}
{Let us prove the sufficiency.} We put
\begin{equation}\label{eq.lemFphi.2}
\Phi_j=\varphi_{\b,\,\d}\cdot h_j\,,\quad j=1,2,
\end{equation}
where $\varphi_{\b,\,\d}$ and $h_j$ are defined by \eqref{eqPhi} and
\eqref{eqh0h1h2}, respectively.

We first show that $\widehat{\Phi}_1\in L_p(\R^d)$. To see this we
use the representation:
\begin{equation}\label{eq.lemFphi.3}
\Phi_1=\Phi_{1,\,1}+\Phi_{1,\,2},
\end{equation}
where
\begin{equation}\label{eq.lemFphi.4}
\Phi_{1,\,1}(x)=\sum_{\nu=1}^{\s-1}\binom{\d}{\nu}(-1)^\nu|x|^{\b
\nu}h_1(x),
\end{equation}
\begin{equation}\label{eq.lemFphi.5}
\Phi_{1,\,2}(x)=h_1(x)\bigg(1+\sum_{\nu=\s}^\infty\binom{\d}{\nu}(-1)^\nu|x|^{\b
\nu}\bigg),
\end{equation}
and $\s> 2(d/p+1)/\b+2$.

From Lemma~\ref{lemRuOdf1} for $\b\not\in 2\N$ and
Lemma~\ref{lemDer} for $\b\in 2\N$ it follows immediately that
$\widehat{\Phi}_{1,\,1}\in L_p(\R^d)$. Using Lemma~\ref{lemDer}
and~\eqref{eqBinom} it is easy to verify that for any positive $\b$
and $\d$
\begin{equation}\label{eq.lemFphi.6}
\widehat{\Phi}_{1,\,2}\in L_p(\R^d).
\end{equation}
Thus, taking into account~\eqref{eq.lemFphi.3}, we obtain that
\begin{equation}\label{eq.lemFphi.7}
\widehat{\Phi}_{1}\in L_p(\R^d).
\end{equation}

Now we check that $\widehat{\Phi}_{2}\in L_p(\R^d)$.  Observe that
for any positive $\b$ and $\d$ the following expansion holds:
\begin{equation}\label{eqRazphi}
(1-|x|^{\b})^\d_+=\sum_{\nu=0}^\infty a_\nu(1-|x|^2)_+^{\d+\nu},
\end{equation}
where $a_\nu\in\R$ and $a_0=(\b/2)^\d$. Consequently, the function
$\Phi_2$ can be represented as follows:
\begin{equation*}\label{eq.lemFphi.8}
\Phi_2=\Phi_{2,\,1}+\Phi_{2,\,2},
\end{equation*}
where
$$
\Phi_{2,\,1}(x)=h_2(x)\sum_{\nu=0}^\l a_\nu (1-|x|^2)_+^{\d+\nu},
$$
$$
\Phi_{2,\,2}(x)=h_2(x)\sum_{\nu=\l+1}^\infty a_\nu
(1-|x|^2)_+^{\d+\nu},
$$
and $\l> d/p-\d$. Using equality~\eqref{eqBesrav} and asymptotic
formula~\eqref{eqBesAs}, it is easy to see that
$\widehat{\Phi}_{2,\,1}\in L_p(\R^d)$. Calculating the  partial
derivatives of the function $\Phi_{2,\,2}$, we obtain that
$\Phi_{2,\,2}\in {\overset{\circ}{W}}\!\!\phantom{l}_1^m(\R^d)$,
where $m=[d/p]+1$. Thus, by Lemma~\ref{lemDer}, we have that
${\widehat{\Phi}}_{2,\,2}\in L_p(\R^d)$. Therefore,
\begin{equation}\label{eq.lemFphi.9}
\widehat{\Phi}_{2}\in L_p(\R^d).
\end{equation}

Combining \eqref{eq.lemFphi.7} and \eqref{eq.lemFphi.9}, we obtain
$\widehat{\varphi}_{\b,\,\d}\in L_p(\R^d)$.

{Now, let us prove the necessity.} From Lemma~\ref{lemProdF} it
follows that
\begin{equation}\label{eq.lemFphi.10}
\widehat{\Phi}_{j}\in L_p(\R^d),\quad j=1,2.
\end{equation}
We claim that under condition \eqref{eq.lemFphi.10}, the pair
$(1/p,\b)$ belongs to $\mathrm{B}(d)$. Indeed, from
\eqref{eq.lemFphi.6}, \eqref{eq.lemFphi.10} and \eqref{eq.lemFphi.3}
it follows immediately that $\widehat{\Phi}_{1,\,1}\in L_p(\R^d)$.
Taking into account that
\begin{equation*}
\Phi_{1,\,1}(x)=|x|^\b\phi(x),
\end{equation*}
where
$$
\phi(x)=\sum_{\nu=1}^{\s-1}\binom{\d}{\nu}(-1)^\nu|x|^{\b
(\nu-1)}h_1(x)
$$
and applying Lemma~\ref{lemRuOdf1} we conclude that
$(1/p,\b)\in\mathrm{B}(d)$.

Now we show that $\d>d(1/p-1/2)-1/2$. Similarly to the previous
arguments, we have $\widehat{\Phi}_{2,\,1}\in L_p(\R^d)$. Using
Lemma~\ref{lemPsiBes}, we get
\begin{equation*}
\begin{split}
\int_1^N|\widehat{\Phi}_{2,\,1}(x)|^{{p}} dx\ge C\int_1^N
r^{d-1}\bigg|&\frac{\cos(2\pi r-\frac{\pi\d}{2}-\frac\pi
4)}{r^{\frac{d+1}{2}+\d}}\bigg|^{{p}} dr-\\
&-O\bigg(\int_1^N r^{d-1-{p}(\frac{d+3}{2}+\d)}dr\bigg).
\end{split}
\end{equation*}
The last inequality implies that $\d>d(1/p-1/2)-1/2$. Otherwise, we
would have that $\widehat{\Phi}_2\not\in L_p(\R^d)$.
%
\end{proof}

%
%

Let us consider general approximation methods generated by the
kernel
\begin{equation*}
\mathcal{K}_n^{\varphi}(x)=\sum_{k}\varphi\bigg(\frac k n\bigg)
e^{i(k,x)},
\end{equation*}
where $\varphi\in C(\R^d)$ is a real valued centrally symmetric
function with a compact support in $\{x\in\R^d\,:\,|x|\le 1\}$ and
$\varphi(0)=1$.  By analogy with the definition of the methods
\eqref{eqFMBR} and \eqref{eqLPBR}, we put:
\begin{equation}\label{eqFMBRF}
\mathcal{L}_n^{\varphi}(f;x)=(2\pi)^{-d}\int_{\T^d}f(x+y)\mathcal{K}_n^{\varphi}(y)dy,\quad
n\in\N,
\end{equation}
\begin{equation}\label{eqLPBRF}
\mathcal{L}_{n;\,\l}^{\varphi}(f;x)=(2n+1)^{-d}\sum_{k=0}^{2n}f(t_n^k+\l)\mathcal{K}_n^{\varphi}(x-t_n^k-\l),\quad
n\in\N.
\end{equation}

As usual, the norm of a linear operator $\mathcal{L}_n^{\varphi}$ is
given by
$$
\Vert \mathcal{L}_n^{\varphi}\Vert_{(p)}=\sup_{\Vert f\Vert_p\le
1}\Vert \mathcal{L}_n^{\varphi}(f)\Vert_p\,.
$$
By analogy, we define the (quasi-)norm of a family
$\{\mathcal{L}_{n;\,\lambda}^{\varphi}\}$ by
$$
\Vert
\{\mathcal{L}_{n;\,\lambda}^{\varphi}\}\Vert_{(p)}=(2\pi)^{-d/p}\sup_{\Vert
f\Vert_p\le 1}\Vert
\mathcal{L}_{n;\,\lambda}^{\varphi}(f;x)\Vert_{\overline{p}}\,.
$$

The proof of the lemma below is standard (see, for
example,~\cite{RS08}).

\begin{lemma}\label{lemSh1}
\emph{1) The means $\mathcal{L}_n^{\varphi}$ converge in  $L_p$,
$1\le p\le\infty$, if and only if the sequence of their norms
$\{\Vert \mathcal{L}_n^{\varphi}\Vert_{(p)}\}_{n\in\N}$ is bounded.}

\emph{2)  The family $\{\mathcal{L}_{n;\,\l}^{\varphi}\}$ converges
in $L_p$, $0<p\le\infty$, if and only if the sequence $\{\Vert
\{\mathcal{L}_{n;\,\lambda}^{\varphi}\}\Vert_{(p)}\}_{n\in\N}$ is
bounded.}
\end{lemma}

The general conditions of convergence for the methods
\eqref{eqFMBRF} and \eqref{eqLPBRF} are formulated in the next lemma
(see~\cite{RS04} and \cite{RRS09}).

\begin{lemma}\label{lemShFLp}
\emph{1) The means $\mathcal{L}_n^{\varphi}$ and the family
$\{\mathcal{L}_{n;\,\l}^{\varphi}\}$ converge in $L_p$ for all $1\le
p\le\infty$ if and only if $\widehat{\varphi}\in L_1(\R^d)$.}

\emph{2) Let $0<p\le 1$. The family
$\{\mathcal{L}_{n;\,\l}^{\varphi}\}$ converges in $L_p$ if and only
if $\widehat{\varphi}\in L_p(\R^d)$.}
\end{lemma}


The following two lemmas are the main tools for proving the theorems
from Section~\ref{s2.2}.

\begin{lemma}\label{lemPolKf}
\emph{Let $\d>0$, $(1/p,\b)\in \mathrm{B}(d)$, and $(
1/p,\delta)\in\Sigma(d)\cup\Omega'(d)$. Then
\begin{equation}\label{eq.PolKf.00}
\Vert T-\mathcal{S}_n^{\b,\,\d}(T)\Vert_p\asymp
n^{-\b}\Vert\Delta^{\b/2}T\Vert_p\,,\quad T\in\mathcal{T}_n,\quad
n\in \N.
\end{equation}
In \eqref{eq.PolKf.00} the operator $\mathcal{S}_n^{\b,\,\d}$ can be
replaced by $\mathcal{S}_{n;\,\l}^{\b,\,\d}$ for any fixed
$\lambda\in\R^d$ without affecting the constants.}
\end{lemma}

\begin{proof}
Let us prove the upper estimate. It is easy to check that for each
polynomial $T\in\mathcal{T}_n$, $n\in\N$ and $\lambda\in\R^d$
\begin{equation}\label{eq.PolKf.1}
\mathcal{S}_n^{\b,\,\d}(T;x)=\mathcal{S}_{n;\,\lambda}^{\b,\,\d}(T;x),\quad
x\in\T^d.
\end{equation}
From \eqref{eq.PolKf.1} it follows that  the operator
$\mathcal{S}_n^{\b,\,\d}$ in~\eqref{eq.PolKf.00} can be replaced by
$\mathcal{S}_{n;\,\l}^{\b,\,\d}$.

According to Lemma~\ref{lemSh1} for $p\ge 1$, we have
\begin{equation}\label{eq.PolKf.2}
\Vert \mathcal{S}_n^{\b,\,\d}(T)\Vert_p\le C\Vert T\Vert_p\,,\quad
T\in\mathcal{T}_n,\quad n\in\N.
\end{equation}
We claim that the inequality \eqref{eq.PolKf.2} holds also for
$0<p<1$. Indeed, from the equality~\eqref{eq.PolKf.1} and the
conditions of Lemma~\ref{lemPolKf} we obtain that for every
polynomial $T\in\mathcal{T}_n$ and $n\in\N$
$$
\Vert \mathcal{S}_n^{\b,\,\d}(T)\Vert_p=(2\pi)^{-d/p}\Vert
\mathcal{S}_{n;\,\lambda}^{\b,\,\d}(T)\Vert_{\overline{p}}\le \Vert
\mathcal{S}_{n;\,\lambda}^{\b,\,\d}\Vert_{(p)}\Vert T\Vert_p\le
C\Vert T\Vert_p\,.
$$
Thus, we have for $0<p\le\infty$
\begin{equation}\label{eq.PolKf.2(1)}
\varphi_{\b,\,\d}\in M_p(\mathcal{T}).
\end{equation}

Put
\begin{equation}\label{eq.PolKf.3}
\xi(x)=\left\{
         \begin{array}{ll}
           {(1-\varphi_{\b,\,\d}(x))}{|x|^{-\b}}, & \hbox{$x\neq 0$;} \\
           \d, & \hbox{$x=0$,}
         \end{array}
       \right.
\end{equation}
$$
\xi_j(x)=h_j(x)\xi(x),\quad j=1,2,
$$
where $h_1$ and $h_2$ are defined by~\eqref{eqh0h1h2}. Note that
$$
\xi_1(x)=h_1(x)\sum_{\nu=1}^\infty\binom{\d}{\nu}(-1)^\nu|x|^{\b
(\nu-1)}.
$$
By analogy with the proof of sufficiency in Lemma~\ref{lemFphi} we
have that $\widehat{\xi}_1\in L_{p^*}(\R^d)$ $(p^*=\min(p,1))$ and
hence, by Lemma~\ref{lemDostNer}, we get
\begin{equation}\label{eq.PolKf.4}
\xi_1\in M_p(\mathcal{T}).
\end{equation}
Now we consider the function $\xi_2$. It is obvious that
\begin{equation}\label{eq.PolKf.55}
(h_2(x)|x|^{-\b})\,\,\widehat{}\,\,\in L_{p^*}(\R^d).
\end{equation}
Applying Lemma~\ref{lemDostNer}, \eqref{eq.PolKf.55},
Lemma~\ref{lemSvMp}, and~\eqref{eq.PolKf.2(1)}, we get
\begin{equation}\label{eq.PolKf.6}
\xi_2\in M_p(\mathcal{T}).
\end{equation}
Thus, by \eqref{eq.PolKf.4}, \eqref{eq.PolKf.6} and
Lemma~\ref{lemSvMp}, we obtain that  $\xi\in M_p(\mathcal{T})$ and
hence the application of Lemma~\ref{lemPSt} yields
\begin{equation}\label{eq.PolKf.7}
\Vert T-\mathcal{S}_n^{\b,\,\d}(T)\Vert_p\le C_3
n^{-\b}\Vert\Delta^{\b/2}T\Vert_p\,,\quad T\in\mathcal{T}_n,\quad
n\in\N.
\end{equation}

In order to prove the lower estimate in~\eqref{eq.PolKf.00} we put
$$
\eta_j(x)=h_j(x)(\xi(x))^{-1},\quad  j=1,2.
$$
We claim that $\eta_j\in M_p(\mathcal{T})$, $j=1,2$. To see that
$\eta_2\in M_p(\mathcal{T})$, we represent $\eta_2$ in the following
form:
$$
\eta_2=\eta_{2,\,1}+\eta_{2,\,2},
$$
where
\begin{equation*}
\eta_{2,\,1}(x)=|x|^{\b}h_2(x)\sum_{\nu=0}^{\sigma}\varphi^\nu(x),\quad
\eta_{2,\,2}(x)=|x|^{\b} h_2(x)\sum_{\nu=\sigma+1}^\infty
\varphi^\nu(x),
\end{equation*}
and $\sigma>(d/p+1)/\delta$. Using Lemmas~\ref{lemDer}
and~\ref{lemDostNer}, it is easy to check that
\begin{equation}\label{eq.PolKf.8}
\eta_{2,\,2}\in M_p(\mathcal{T})
\end{equation}

Consider the function $\eta_{2,\,1}$. It is obvious that
\begin{equation}\label{eq.PolKf.9}
(h_2(x)|x|^{\b})\,\,\widehat{}\,\,\in L_{p^*}(\R^d).
\end{equation}
Applying Lemma~\ref{lemDostNer}, \eqref{eq.PolKf.9},
Lemma~\ref{lemSvMp}, and~\eqref{eq.PolKf.2(1)}, we get
\begin{equation}\label{eq.PolKf.10}
\eta_{2,\,1}\in M_p(\mathcal{T}).
\end{equation}
Thus, by \eqref{eq.PolKf.8}, \eqref{eq.PolKf.10} and
Lemma~\ref{lemSvMp}, we obtain
\begin{equation}\label{eq.PolKf.10}
\eta_{2}\in M_p(\mathcal{T}).
\end{equation}

Now to check
\begin{equation}\label{eqEt11Mp}
\eta_{1}\in M_p(\mathcal{T}),
\end{equation}
we introduce the function
\begin{equation}\label{PhiF}
\phi(x)=\frac{|x|^{\b}}{1-(1-|x|^\b)_+^\d}-\sum_{\nu=0}^{\s}a_{\nu}|x|^{\nu\b},
\end{equation}
where $\s>\frac 2\b(\frac dp+1)+2$. We put also $
\gamma(x)=\phi(x)h_1(x). $ From Lemmas~\ref{lemSvMp}
and~\ref{lemRuOdf1} it  follows that ${\gamma}-{\eta}_{1}\in
M_p(\mathcal{T})$. Thus, to check~\eqref{eqEt11Mp}, we need to check
\begin{equation}\label{eqG}
{\gamma}\in M_p(\mathcal{T}).
\end{equation}

Note that, for $u\in(0,1)$
$$
\frac{u^\b}{1-(1-u^\b)^\d}=\bigg(\sum_{\nu=0}^\infty c_\nu u^{\b
\nu}\bigg)^{-1},
$$
where $c_0=\d\neq 0$. Therefore, the numbers $\{a_\nu\}$
in~\eqref{PhiF} can be chosen such that
\begin{equation}\label{eqth4vsp1}
\phi(x)={\bigg(\sum_{\mu=0}^\infty c_{\mu} |x|^{\b
\mu}\bigg)^{-1}}{\sum_{\nu=\s+1}^\infty {b}_{\nu}|x|^{\b \nu}}.
\end{equation}
Calculating the partial derivatives of $\phi(x)$ and taking into
account \eqref{eqth4vsp1}, we get $h_1\phi\in
{\overset{\circ}{W}}\!\!\phantom{l}_1^m(\R^d)$, where $m=[d/p]+1$.
Thus, by Lemma~\ref{lemDer} we have \eqref{eqG} and therefore
\eqref{eqEt11Mp}.

Combining \eqref{eq.PolKf.10} and \eqref{eqEt11Mp}, using
Lemmas~\ref{lemSvMp} and~\ref{lemPSt}, we have the following
inequality:
\begin{equation}\label{eq.PolKf.7(2)}
n^{-\b}\Vert\Delta^{\b/2}T\Vert_p \le C_4\Vert
T-\mathcal{S}_n^{\b,\,\d}(T)\Vert_p\,,\quad T\in\mathcal{T}_n,\quad
n\in\N.
\end{equation}
Thus, from \eqref{eq.PolKf.7} and \eqref{eq.PolKf.7(2)} we have the
two-sided inequality~\eqref{eq.PolKf.00}.
%
\end{proof}

\begin{lemma}\label{lemPolKf2}
\emph{Let $f\in L_p(\T^d)$, $1\le p\le \infty$, $\b,\d>0$,  and $(
1/p,\delta)\in\Sigma(d)\cup\Omega'(d)$. Then
\begin{equation*}\label{eq.PolKf.0}
\Vert f-\mathcal{S}_n^{\b,\,\d}(f)\Vert_p\le C
n^{-\b}\Vert\Delta^{\b/2}f\Vert_p\,,\quad n\in \N,
\end{equation*}
where $C$ is a constant independent of  $f$ and $n$.}
\end{lemma}

\begin{proof}
The proof of Lemma~\ref{lemPolKf2} is similar to the proof of
Lemma~\ref{lemPolKf} (the upper inequality in~\eqref{eq.PolKf.00}).
The main difference is that it is necessary to use multipliers of
Fourier series instead of multipliers of trigonometric polynomials.
We give a brief proof of Lemma~\ref{lemPolKf2} by using the
corresponding theorems in~\cite{Tr80} (see also Chapters 7 and 8
in~\cite{TrigBook}).


Denote by $M_p$ the algebra of locally Riemann-integrable functions
with the norm
$$
\Vert \varphi \Vert_{M_p}=\sup_n \bigg\Vert\left\{\varphi
\left(\frac kn\right)\right\} \bigg\Vert_{L_p\mapsto L_p}<\infty,
$$
where $\Vert \{\l_k\}\Vert_{L_p\mapsto L_p}$ is the norm of Fourier
multiplier $\Lambda=\{\l_k\}$, acting from $L_p$ to $L_p$ (for the
precise definition see~\cite{Tr80}). By the comparison principle for
Fourier multipliers (see Theorem~6 in~\cite{Tr80}) it suffices to
show that the function $\xi$ defined by~\eqref{eq.PolKf.3} belongs
to $M_p$.

Put $\xi_1(x)=h_1(x)\xi(x)$ and $\xi_2(x)=(1-h_1(x))\xi(x)$. In the
proof of Lemma~\ref{lemPolKf} it was shown that $\widehat{\xi}_1\in
L_1(\R^d)$, hence using Theorem~1 in~\cite{Tr80} and the inequality
$$
\Vert \{\l_k\}\Vert_{L_p\mapsto L_p}\le \Vert
\{\l_k\}\Vert_{L_\infty\mapsto L_\infty},
$$
which holds for any $p\in [1,\infty)$ (see, for example,~\cite[p.
284]{Zig}), we get $\xi_1\in M_p$.

To conclude the proof, it remains to show that $\xi_2\in M_p$. In
accordance with the conditions of Lemma~\ref{lemPolKf2}, we obtain
$\varphi_{\b,\,\d}\in M_p$. Note also that the function
$\psi(x)=(1-h_1(x))|x|^{-\b}$ can be represented as an absolutely
convergent Fourier integral (see Theorem~4 in~\cite{Tr80}), hence
$\psi\in M_p$ (see Theorem~1 in~\cite{Tr80}). Thus, by using
elementary properties of multipliers, we obtain $\xi_2\in M_p$.
%
\end{proof}

\section{Proofs of the Main Results}\label{s4}

\begin{proof}[Proof of Theorem~\ref{thSh1}]
Let $\a$ and $\d$ be some positive numbers. Let us show that the
convergence of the means $\mathcal{S}_{n}^{\a,\,\d}$ in $L_p$
implies the convergence of $\mathcal{S}_{n}^{\b,\,\d}$ in the same
space $L_p$ for any $\b>0$.

First we show that the convergence of $\mathcal{S}_{n}^{\a,\,\d}$
implies the convergence of $\mathcal{S}_{n}^{\a,\,\d+1}$. Indeed,
from the equality
$$
\mathcal{S}_{n}^{\a,\,\d+1}=
\mathcal{S}_{n}^{\a,\,\d}-n^{-\a}\Delta^{\a/2}\circ\mathcal{S}_{n}^{\a,\,\d}
$$
and from Lemma~\ref{lemRuOdfBer}, we get that for each $f\in L_p$\,:
\begin{equation}\label{eq.pr.th1.1}
\begin{split}
\Vert \mathcal{S}_{n}^{\a,\,\d+1}(f)\Vert_p&\le \Vert
\mathcal{S}_{n}^{\a,\,\d}(f)\Vert_p+n^{-\a}\Vert
\Delta^{\a/2}\,\mathcal{S}_{n}^{\a,\,\d}(f)\Vert_p\le \\
&\le C\Vert \mathcal{S}_{n}^{\a,\,\d}(f)\Vert_p\,.
\end{split}
\end{equation}
From the last inequality and Lemma~\ref{lemSh1} it follows that the
convergence of $\mathcal{S}_{n}^{\a,\,\d}$ implies the convergence
of $\mathcal{S}_{n}^{\a,\,\d+1}$.

Next, using the expansion
$$
(1-|x|^\b)_+^\d=\sum_{\nu=0}^\infty b_\nu (1-|x|^\a)_+^{\d+\nu},
$$
we can represent the means $\mathcal{S}_{n}^{\b,\,\d}$ in the
following form
\begin{equation}\label{eq.pr.th1.2}
\mathcal{S}_{n}^{\b,\,\d}=\sum_{\nu=0}^{\l} b_\nu
\mathcal{S}_{n}^{\a,\,\d+\nu}+\mathcal{P}_n,
\end{equation}
where
$$
\mathcal{P}_n(f;x)=\sum_k \psi\left(\frac kn\right)c_k(f)e^{i(k,x)}
$$
and
\begin{equation}\label{eqForLam}
\psi(x)=(1-|x|^\b)_+^\d-\sum_{\nu=0}^{\l}b_\nu
(1-|x|^\a)_+^{\d+\nu}.
\end{equation}
We choose the parameter $\l$ in \eqref{eqForLam} such that
$\psi(0)\neq 0$ and $\l>d$.  Repeating the proof of
Lemma~\ref{lemFphi} it is easy to check that $\widehat{\psi}\in
L_1(\R^d)$. Whence, by Lemma~\ref{lemShFLp}, we have that
$\mathcal{P}_n$ with an appropriate normalization converge in $L_q$
for any $q\in [1,\infty]$. Thus, by using the
equality~\eqref{eq.pr.th1.2}, the inequality~\eqref{eq.pr.th1.1},
and Lemma~\ref{lemSh1}, we get that the means
$\mathcal{S}_{n}^{\b,\,\d}$ converge in $L_p$ for any $\b>0$. To
complete the proof of the theorem, it remains only to use
Lemma~\ref{lemBR2}.
%
\end{proof}

\bigskip


\begin{proof}[Proof of Theorem~\ref{thSh2}]
The proof of item 1) of Theorem~\ref{thSh2} is similar to the proof
of Theorem~\ref{thSh1}. The proof of item 2) follows from
Lemmas~\ref{lemFphi} and~\ref{lemShFLp}.
%
\end{proof}

\bigskip


\begin{proof}[Proof of Theorem~\ref{thKf1}]
By using Lemma~\ref{lemPolKf2} for any $g$ satisfying
$\Delta^{\b/2}g\in L_p(\T^d)$, we get
\begin{equation}\label{eq.Kf1.1}
\begin{split}
\Vert f-\mathcal{S}_n^{\b,\,\d}(f)\Vert_{p}&\le
\Vert (f-g)-\mathcal{S}_n^{\b,\,\d}(f-g)\Vert_{p}+\Vert g-\mathcal{S}_n^{\b,\,\d}(g)\Vert_{p}\le\\
&\le(1+\Vert \mathcal{S}_n^{\b,\,\d}\Vert_{(p)})\Vert
f-g\Vert_{p}+C_1 n^{-\b}\Vert \Delta^{\b/2} g\Vert_{p}\,.
\end{split}
\end{equation}
Passing to the infimum on $g$ in~\eqref{eq.Kf1.1}, we obtain
\begin{equation}\label{eq.Kf1.2}
\Vert f-\mathcal{S}_n^{\b,\,\d}(f)\Vert_{p} \le C_2 K_\b(f,1/n)_p\le
C_2 \widetilde{K}_\b(f,1/n)_p\,.
\end{equation}

We now prove the lower estimate. Using the lower estimate
in~\eqref{eq.PolKf.00}, we have
\begin{equation}\label{eq.Kf1.3}
\begin{split}
K_{\b}(f,1/n)_p &\le\widetilde{K}_{\b}(f,1/n)_p\le
\\
&\le\Vert f-\mathcal{S}_n^{\b,\,\d}(f)\Vert_p+n^{-\b}\Vert
\Delta^{\b/2} \mathcal{S}_n^{\b,\,\d}(f)\Vert_p\le
\\
&\le \Vert f-\mathcal{S}_n^{\b,\,\d}(f)\Vert_p+C_3\Vert
\mathcal{S}_n^{\b,\,\d}(f-\mathcal{S}_n^{\b,\,\d}(f))\Vert_p\le\\
&\le C_4\Vert f-\mathcal{S}_n^{\b,\,\d}(f)\Vert_p\,.
\end{split}
\end{equation}
The equivalence
\begin{equation}\label{eq.Kf1.4}
\Vert f-\mathcal{S}_{n;\,\l}^{\b,\,\d}(f)\Vert_{\overline{p}}\asymp
\Vert f-\mathcal{S}_{n}^{\b,\,\d}(f)\Vert_{{p}}
\end{equation}
follows directly from Lemma~2.2 in~\cite{RRS09} (see also the proof
of Theorem 2 in~\cite{RS08}).

Thus, the equivalences in \eqref{eq.thKf1.0} follow from
\eqref{eq.Kf1.2}, \eqref{eq.Kf1.3}, \eqref{eq.Kf1.4} and Theorem~A.
%
\end{proof}

\bigskip


\begin{proof}[Proof of Theorem~\ref{thKf2}]
The proof of Theorem~\ref{thKf2} follows from Lemma~\ref{lemPolKf}
and Theorem~\ref{thSh2}. It coincides with the proof of Theorem~3 in
\cite{RS08}.
\end{proof}

\bigskip


\begin{proof}[Proof of Theorem~\ref{thOm1}]
The proof of inequality~\eqref{eqOm1} is similar to the proof of
Theorem~6 in~\cite{RS08}.
\end{proof}

\bigskip

\begin{proof}[Proof of Theorem~\ref{thOm2} and Theorem~\ref{thOm3}]
The proof is similar to the proof of Theorem~\ref{thKf1}.
\end{proof}

\end{document}